\documentclass[10pt]{amsart}
\usepackage{amsfonts}
\usepackage{amsthm}
\usepackage{amsmath}
\usepackage{amssymb}
\usepackage{amscd}
\usepackage{float}
\usepackage{mathtools}
\usepackage[centerlast,small]{caption}
\usepackage{longtable}
\usepackage{bm}
\usepackage[latin2]{inputenc}
\usepackage{t1enc}
\usepackage[mathscr]{eucal}
\usepackage{indentfirst}
\usepackage{graphicx}
\usepackage{graphics}
\usepackage{pict2e}
\usepackage{epic}
\usepackage[margin=2.9cm]{geometry}
\usepackage{tikz}
\usepackage[hidelinks]{hyperref}
\usepackage{boldline}
\usepackage{makecell}
\usepackage{filecontents}
\usepackage[numbers]{natbib}

\theoremstyle{definition}
\newtheorem{Th}{Theorem}[section]

\newtheorem{Prop}[Th]{Proposition}
\newtheorem{Def}[Th]{Definition}

\newtheorem{Rem}[Th]{Remark}
\newtheorem{?}[Th]{Problem}
\newtheorem{Lem}[Th]{Lemma}

\newcommand{\ovl}{\overline}

\newcommand{\Z}{\mathbb{Z}}

\newcommand{\tb}{\textbf}
\newcommand{\ti}{\textit}
\newcommand{\mcl}{\mathcal}
\newcommand{\mfk}{\mathfrak}
\newcommand{\s}{\sigma}
\newcommand{\be}{\begin{equation*}}
\newcommand{\ee}{\end{equation*}}

\title {Knots with Prism Manifold Surgeries}
\author{Zhengyuan Shang}
\date{\today} 

\begin{document}
\maketitle

\begin{abstract}
Ballinger et al. have determined the list of all prism manifolds that are possibly realizable by Dehn surgeries on knots in $S^3$. In this paper, we explicitly find braid words of primitive/Seifert-fibered knots on which surface slope surgeries yield all the prism manifolds listed above. This completes the solution to the prism manifold realization problem.
\end{abstract}

\section{Introduction}
Every closed orientable three manifold admits a unique decomposition into prime three manifolds, each of which can be uniquely assigned one of eight possible geometric structures by the Geometrization theorem. The simplest three manifolds are the ones with finite fundamental groups, or equivalently, the spherical manifolds. Since all compact orientable three manifolds can be realized by Dehn surgeries on links in $S^3$ \cite{WL, AW}, it is of interest to study which spherical manifolds are realizable by surgeries on knots.\\

A knot in $S^3$ is either torus, satellite, or hyperbolic. Finite surgeries on torus knots are classified by Moser \cite{LM}. Then Bleiler and Hodgson classified finite surgeries on satellite knots \cite{BH}, while Boyer and Zhang showed that the coefficient of any finite surgery is either integral or half integral \cite{BZ}. Later, Li and Ni showed every spherical manifold resulting from half integral surgery on a hyperbolic knot can also be obtained through surgery on a torus knot or a satellite knot \cite{LN}. Therefore, we restrict our attention to integral surgeries, or equivalently, positive integral surgeries.   \\

Spherical three manifolds are classified into: cyclic($\tb{C}$), dihedral($\tb{D}$), tetrahedral($\tb{T}$), octahedral($\tb{O}$), and icosahedral($\tb{I}$) types. The $\tb{C}$-type manifolds realization problem is addressed by Greene \cite{JG}, then Gu \cite{LG} solves the cases for $\tb{T}$, $\tb{O}$, $\tb{I}$-type manifolds, which leaves $\tb{D}$-type manifolds, also known as prism manifolds, the only remaining case. We denote as $\mcl{P}(p, q)$ the prism manifold with Seifert invariants
\be
\label{eq:def}
(-1; (2, 1), (2, 1), (p, q))
\ee
where $p>1$ and gcd$(p, q)=1$. Then by \cite[Lemma 6.2]{NZ}, $\mcl{P}(p, q)$ could only be realized by $4|q|$ surgeries on knots. Following an approach similar to Greene's, Ballinger et al. have made a list of all possibly realizable prism manifolds $\mcl{P}(p, q)$, with conditions that $q<0, q>p,$ and $0<q<p$ in \cite{B1}, \cite{B2}, \cite{B3} respectively. It is expected that their list is exact. Our main result confirms this conjecture:

\Th{A prism manifold $\mcl{P}(p, q)$ is realizable by a positive integral surgery on a knot in $S^3$ if and only if it appears in \autoref{tab:list}. In particular, $4|q|$-surgeries on knots specified by the listed braid words yield the corresponding prism manifolds.}\\

\tb{Notations} $K$ is always a knot of the form $(\s_1...\s_{B-1})^k(\s_1...\s_T)^a$ or $(\s_1...\s_{B-1})^k(\s_T...\s_1)^a$, where $0<T< B-1$, $k>0$, and $a\in \Z$. The braid word  is considered as part of the structure of the knot.

\begin{figure}[H]
\centering
\includegraphics[scale=0.4]{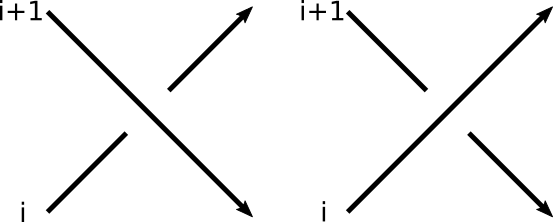}
\caption{Conventions for the braid group: $\sigma_i$ (left) and $\s_i^{-1}$ (right)}
\end{figure}

 Let $\mcl{T}(m, n)$ be the torus link that winds $m$ times around the longitude and $n$ times around the meridian. Then $(\s_1...\s_{B-1})^k$ is $\mcl{T}(B, k)$ and $(\s_1...\s_T)^a$ or $(\s_T...\s_1)^a$ is $\mcl{T}(T+1, a)$. By cutting open the first $T+1$ parallel strands along the longitudes of $\mcl{T}(B, k)$ and $\mcl{T}(T+1, a)$ respectively and gluing the torus tangles together with appropriate choices of orientations, we obtain an embedding of $K$ in the genus two surface $\Sigma$. This will be our embedding for knots in families Spor, $1B_1$, $2$, $3B$, 4, and 5. Let the handlebodies inside and outside $\Sigma$ be $H$ and $H'$ respectively. Suppose $\pi_1(H)$ is generated by $x, y$ and $\pi_1(H')$ is generated by $x', y'$. Let $w, w'$ be the conjugacy classes of $K$ in $\pi_1(H)$ and $\pi_1(H')$. Note that the exponents of each generator in $w$ or $w'$ always have the same sign. We choose $x, y, x', y'$ such that their exponents in $w, w'$ are all positive. 

\begin{figure}[H]
\centering
\includegraphics[scale=0.4]{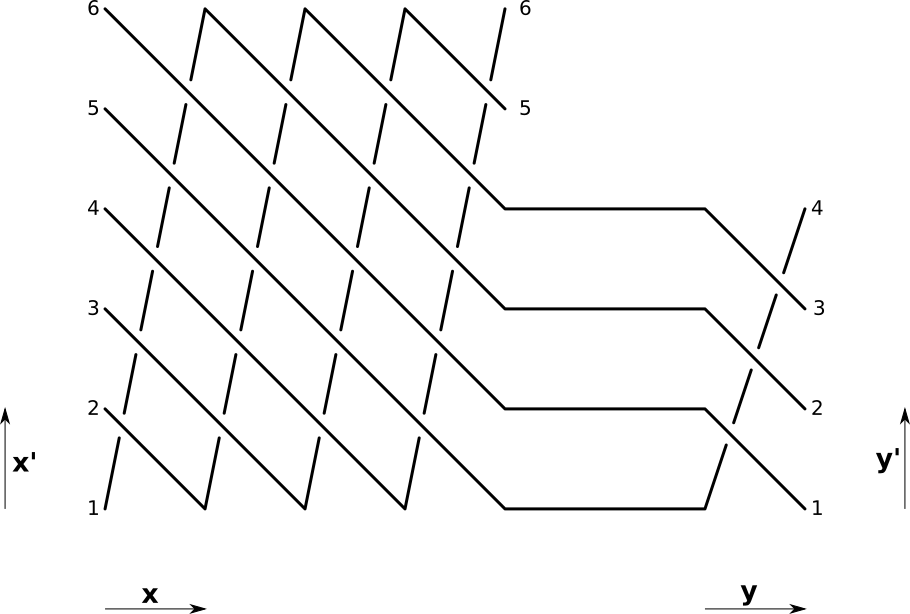}
\caption{The knot $(\s_1...\s_5)^4\s_1...\s_3$  of family $5^+_1$}
\label{fig:fig2}
\end{figure}

For example, the knot above is obtained by gluing together $\mcl{T}(6, 4)$ and $\mcl{T}(4, 1)$. Here $x, x'$ are the longitudes of these torus links, so they generate $\pi_1(H)$. We have $w=(xy)^2(x^2y)^2$ and $w'=x'^4y'$. \\

For families $1A$, $1B_2$, and $3A$, we cut open the first $T+1$ meridional strands and the $T+1$ longitudinal strands of $\mcl{T}(B, k)$ and $\mcl{T}(T+1, a)$ respectively and glue the tangles together to obtain an embedding of $K$ in $\Sigma$. Suppose $\pi_1(H)$, $\pi_1(H')$ are still generated by $x, y$ and $x', y'$. Then this new embedding simply amounts to switching the position of $x$ and $x'$ (cf. \autoref{fig:fig3}).\\

\tb{Organization} We first make some preliminary observations about knots in Table 1 and address families Spor, $1A$, and $5_\ast$ in Section 2. Then we prove the primitive/Seifert-fibered properties of all other families in Section 3 and conclude that $4|q|$ surgeries on them yield prism manifolds $\mcl{P}(p, q)$. In Section 4, we outline the approach to determine $|p|$ for each family.\\

\tb{Acknowledgement} I want to thank my mentor Prof. Yi Ni for his proposal of this topic and many pertinent suggestions. This project would not be possible without his guidance. I am sponsored by the Caltech Summer Undergraduate Research Fellowships (SURF) program and wish to thank Arthur R. Adams for his generous support through SURF. I also want to thank William Ballinger for his program to compute Alexander polynomials from the changemakers \cite[Definition 1.4]{B1} of knots.\\

\section{The Braid Words} 

We combine \cite[Table 5]{B1}, \cite[Table 3]{B2}, and \cite[Theorem 1.1 and 1.2]{B3} to obtain the parameterizations of all realizable prism manifolds in the table below. We keep Ballinger et al.'s original classification of such knots into families $1A, 1B, 2, 3A, 3B, 4, 5$ and Spor. For non-sporadic knots, the superscript identifies the sign of $q$ and the occasional subscript of 1, 2, or $\ast$ clarifies its subfamily. Here, $p$ and $q$ are parametrized by integers $s$ or $s, t$, whose ranges are specified. The parameter $r$ originates from a more concise classification of realizable prism manifolds \cite[Tables 1 and 2]{B1}.\\
\\

\begin{longtable}{ c c c }
\caption{Matching Prism Manifolds with Braid Words}\\

\hlineB{3}
\tb{Prism Manifold} & \tb{Braid Word}& \textbf{$\mathcal{P}$ type} \\
\hlineB{3}
\endfirsthead
\hlineB{3}
\tb{Prism Manifold} & \tb{Braid Word}& \textbf{$\mathcal{P}$ type} \\
\hlineB{3}
\endhead

$p=11, -q=30$ & $(\s_1...\s_{17})^{7}(\s_2\s_1)^{-2}$ & Spor\\
\hline
$p=17, -q=31$ & $(\s_1...\s_9)^{13}(\s_2\s_1)^{-2}$ & Spor\\
\hline 
$p=13, -q=47$ & $(\s_1...\s_{17})^9(\s_{12}...\s_1)^2$ & Spor\\ 
\hline 
$p=23, -q=64$ & $(\s_1...\s_{13})^{19}(\s_4...\s_1)^{-2}$ & Spor\\ 
\hline 
$p=11, q=19$ & $(\s_1...\s_6)^{10}(\s_1\s_2)^2$ &  Spor\\ 
\hline 
$p=13, q=34$ & $(\s_1\s_2...\s_8)^{14}(\s_1...\s_4)^2$ &  Spor\\ 

\hlineB{3}
\makecell{$p=2s+5$\\ $-q=2s^2+7s+7$\\$s\geq 1$} & $(\s_1...\s_{2s+2})^{4s+8}(\s_1)^2$ & $1A^-$\\ 
\hline 
\makecell{$p=2s-1$\\ $q=2s^2+s+1$\\$s\geq 2$} & $(\s_1...\s_{2s})^{4s}(\s_1)^2$ & $1A^+$\\ 
\hline 
\makecell{$p=22s+39$\\ $-q=22s^2+75s+64$\\$s\geq 0$} & $(\s_1...\s_{8s+13})^{12s+20}(\s_1...\s_{4s+7})^{-2s-3}$ & $1B^{-}_1$\\ 
\hline 
\makecell{$p=22s+27$\\ $q=22s^2+57s+37$\\$s\geq 0$} & $(\s_1...\s_{8s+9})^{12s+16}(\s_1...\s_{4s+3})^{-2s-3}$ & $1B^+_1$\\ 
\hline 
\makecell{$p=22s+41$\\ $-q=22s^2+79s+71$\\$s\geq 0$} & $(\s_1...\s_{8s+13})^{12s+22}(\s_{4s+7}...\s_1)^{-2s-3}$ & $1B^-_2$\\ 
\hline 
\makecell{$p=22s+25$\\ $q=22s^2+53s+32$\\$s\geq 0$} & $(\s_1...\s_{8s+9})^{12s+14}(\s_{4s+3}...\s_1)^{-2s-3}$ & $1B^+_2$\\ 
\hlineB{3}

\makecell{$p=16ts+8s+30t+11$\\ $-q=16ts^2+56ts+8s^2+24s+49t+18$\\$r=4s+7, s\geq 0, t\geq 0$} & 
$(\s_1...\s_{8s+13})^{8ts+4s+14t+6}(\s_1...\s_{8s+11})^{-1}$
& $2^-_1$\\ 
\hline 
\makecell{$p=16ts+8s+18t+13$\\ $q=16ts^2+40ts+8s^2+24s+25t+18$\\$r=-4s-5, s\geq 1, t\geq 0$} & 
$(\s_1...\s_{8s+9})^{8ts+4s+10t+8}(\s_1...\s_{8s+7})^{-1}$
& $2^+_1$\\ 
\hline 
\makecell{$p=16ts+8s+18t+5$\\ $-q=16ts^2+40ts+8s^2+16s+25t+7$\\$r=-4s-5, s\geq 0, t\geq 0$} & 
$(\s_1...\s_{8s+9})^{8ts+4s+10t+2}(\s_1...\s_{8s+7})$ & $2^-_2$\\ 
\hline 
\makecell{$p=16ts+8s+30t+19$\\ $q=16ts^2+56ts+8s^2+32s+49t+31$\\$r=4s+7, s\geq 0, t\geq 0$} & 
$(\s_1...\s_{8s+13})^{8ts+4s+14t+8}(\s_1...\s_{8s+11})$
& $2^+_2$\\ 
\hlineB{3}


\makecell{$p=4ts+2s+8t+3$\\ $-q=4ts^2+16ts+2s^2+11s+16t+14$\\$r=2t+1, s\geq 0, t\geq 0$} & 
$(\s_1...\s_{4s+7})^{4ts+2s+8t+5}(\s_{4s+7}...\s_1)^2$
& $3A^-$\\ 
\hline 
\makecell{$p=4ts+6s+8t+13$\\ $q=4ts^2+16ts+6s^2+21s+16t+18$\\$r=2t+3, s\geq 0, t\geq 1$} & 
$(\s_1...\s_{4s+7})^{4ts+6s+8t+11}(\s_{4s+7}...\s_1)^{-2}$ 
& $3A^+$\\ 
\hline 
\makecell{$p=4ts+6s+10t+11$\\ $-q=4ts^2+20ts+6s^2+27s+25t+30$\\$r=2t+3, s\geq -1, t\geq 1$} & $(\s_1...\s_{4s+9})^{4ts+4s+10t+8}(\s_1...\s_{4s+7})^{2s+5}$ & $3B^-$\\ 
\hline 
\makecell{$p=4ts+2s+6t+7$\\ $q=4ts^2+12ts+2s^2+9s+9t+9$\\$r=2t+1, s\geq 0, t\geq 0$} & $(\s_1...\s_{4s+5})^{4ts+4s+6t+8}(\s_1...\s_{4s+3})^{-2s-3}$ & $3B^+$\\ 
\hlineB{3} 

\makecell{$p=8ts^2+8ts+8s+2t+3$\\ $-q=16ts^2+24ts+16s+9t+14$ \\$r=2s+1, s\geq 1 , t\geq 0$} & $(\s_1...\s_{8s+5})^{8ts+6t+8}(\s_1...\s_7)$ & $4^-_1$\\ 
\hline 
\makecell{$p=8ts^2+24ts+8s+18t+13$\\ $q=16ts^2+40ts+16s+25t+18$\\ $r=-2s-3, s\geq 0, t\geq 0$} & $(\s_1...\s_{8s+9})^{8ts+10t+8}(\s_1...\s_7)^{-1}$ & $4^+_1$\\ 
\hline
\makecell{$p=8ts^2+40ts+8s^2+32s+50t+29$\\ $-q=16ts^2+72ts+16s^2+56s+81t+47$\\$r=-2s-5, s\geq 0, t\geq 0$} & $(\s_1...\s_{8s+17})^{8ts+8s+18t+10}(\s_1...\s_7)$ & $4^-_2$\\ 
\hline 
\makecell{$p=8ts^2+24ts+8s^2+16s+18t+7$\\ $q=16ts^2+56ts+16s^2+40s+49t+19$\\ $r=2s+3, s\geq 0, t\geq 0$} & $(\s_1...\s_{8s+13})^{8ts+8s+14t+6}(\s_1...\s_7)^{-1}$ & $4^+_2$\\ 
\hlineB{3}

\makecell{$p=4ts^2+8ts+4s+2t+1$\\ $-q=4ts^2+12ts+4s+9t+5$\\ $r=2s+3,s\geq0, t\geq 0, (s, t)\neq (0, 0)$} & $(\s_1...\s_{4s+5})^{4ts+6t+4}(\s_1...\s_{3})^{-1}$ & $5^-_1$\\
\hline 
\makecell{$p=4s^2t+16st+4s+14t+11$\\ $q=4s^2t+12st+4s+9t+7$ \\$r=-2s-3, s\geq 0, t\geq 0$} & $(\s_1...\s_{4s+5})^{4ts+6t+4}(\s_1...\s_{3})$ & $5^+_1$\\ 
\hline 
\makecell{$p=4ts^2+16ts+4s^2+12s+14t+3$\\ $-q=4ts^2+12ts+4s^2+8s+9t+2$\\ $r=-2s-3,s\geq 0, t\geq 0$} & $(\s_1...\s_{4s+5})^{4ts+4s+6t+2}(\s_1...\s_3)^{-1}$ & $5^-_2$\\ 
\hline 
\makecell{$p=4ts^2+16ts+4s^2+12s+14t+9$\\ $q=4ts^2+20ts+4s^2+16s+25t+16$\\ $r=2s+5, s\geq -1, t\geq 0, (s, t)\neq (-1, 0)$} & $(\s_1...\s_{4s+9})^{4ts+4s+10t+6}(\s_1...\s_{3})$ & $5^+_2$\\ 
\hline
\makecell{$p=-2q-1$ \\$r=-1, -q\geq2$} & $\s_1^{-2q-1}$ & $5^-_*$\\ 
\hline 
\makecell{$p=2q+1$\\$r=-1, q\geq 1$} & $\s_1^{2q+1}$ & $5^+_*$\\
\hlineB{3}
\label{tab:list}
\end{longtable}

\Rem Except for the torus knots in family $5_\ast$, all knots above resemble the twisted torus knots described in \cite{JD}, which provide a rich class of small Seifert-fibered spaces (three exceptional fibers or fewer) surgeries. In fact, our search for the desired knots starts from "generalized" twisted torus knots of the forms $(\s_1...\s_{B})^k(\s_1...\s_{T+1})^a$ or $(\s_1...\s_{B})^k(\s_{T+1}...\s_{1})^a$. For each family of knots above, we first compute their Alexander polynomials through the changemakers given in \cite{B1, B2, B3}. Then with Mathematica, we look for "generalized" twisted torus knots (usually with $B\leq 20$, $k\leq 30$) that have the same Alexander polynomials. By performing surgeries on such knots in small cases ($1\leq s, t\leq 3$) and comparing the resulting manifolds with the expected ones, we are able to obtain a few candidates for each family, which eventually leads to the tabulation above.

\Rem By \cite[Lemma 2.1]{B1}, for $m\geq 1$, $4m$ surgeries on $\mcl{T}(2m+1, 2)$ give the manifolds $\mcl{P}(2m+1, m)$. Hence for $q\geq 1$, $4|q|$ surgeries on $\s_1^{2q+1}$ of family $5_\ast^+$ yield the desired prism manifolds. The same lemma implies for $q\leq -2$, $4|q|$ surgeries on $\s_1^{-2q-1}$ of family $5_\ast^-$ result in the prism manifolds $\mcl{P}(-2q-1, q)$. Moreover, for the six sporadic knots, we have verified with SnapPy that $4|q|$ surgeries on them also give rise to their corresponding prism manifolds.

\Lem{Suppose $\mcl{P}(p, q)$ is a realizable with $p\in \Z$ and $q>0$. (Here for $p<0$, $\mcl{P}(p, q)=\mcl{P}(-p, -q)$ as defined in the introduction). Then $p= 3$ (mod 4) if $q$ is odd and $p= 1$ (mod 4) if $q$ is even.}

\begin{proof}
Let $\Delta_K(\mfk{t}), S_{4q}^3(K), \lambda(\mcl{P}(p, q))$, and $\mathfrak{s}(p, q)$ be the Alexander polynomial of $K$ in $\mfk{t}$, the manifold obtained from $4q$ surgery on $K$, the Casson-Walker invariant of $\mcl{P}(p, q)$, and the Dedekind sum in $p$ and $q$. By \cite[Proposition 6.1.1]{CL} and the reciprocity law of Dedekind sums, we have 
\be
\lambda(\mcl{P}(p, q))=\frac{p}{8q}-\mathfrak{s}(p, q)
\ee
Meanwhile, \cite[Theorem 2.8]{BL} implies 
\be
\lambda(S^3_{4q}(K))=-\frac{(2q-1)(4q-1)}{24q}+\frac{1}{4q}\Delta_K''(1)
\ee
Suppose $S_{4q}^3(K)=\mcl{P}(p, q)$, then we have 
\be
3p-24q\mathfrak{s}(p, q)=-(2q-1)(4q-1)+6\Delta_K''(1) \ \ \ \ (\text{mod } 24q)
\ee
which simplifies to 
\be
1-6q+8q^2+p(-1+6q-2q^2)=6\Delta_K''(1)\ \ \ \  (\text{mod } 24q)
\ee
Since $\Delta_K(\mfk{t})=\Delta_K(-\mfk{t})$, $\Delta_K''(1)$ is even. Hence $4|1-6q+p(-1+6q-2q^2)$. If $q$ is even, we have $4|1-p$, so $p=1$ (mod 4). If $q$ is odd, $4|1-6+p(-1+6-2)=3p-5$, so $p=3$ (mod 4).
\end{proof}

\Rem This is a generalization to \cite[Lemma 2.3]{B3}, and we follow a similar proof approach. This lemma allows us to decide the sign of a realizable prism manifold once we have determined $|p|$ and $|q|$. For all the prism manifolds $\mcl{P}(p, q)$ in Table \ref{tab:list}, we have verified that $p=1$ (mod 4) when $q$ is even and $p= 3$ (mod 4) when $q$ is odd, which means to prove Theorem 1.1, we only need to show that $4|q|$ surgeries on $K$ yield prism manifolds of three exceptional fibers of multiplicities $(2, 2, |p|)$. 

\Lem The surface slopes of knots $K$ in Table 1 (except for family $5_\ast$ knots) are $4|q|$. Hence $4|q|$ surgeries are also surface slope surgeries.

\begin{proof}
Note that the surface slope of the torus link $\mcl{T}(m, n)$ is $mn$. Since the two tangles contribute $Bk$ and $(T+1)a$ respectively, the slope of $K$ with respect to $\Sigma$ is the sum $Bk+(T+1)a$. For all knots in Table 1 (except for $5_\ast$), we verify explicitly that $4|q|=Bk+(T+1)a$. 
\end{proof}

\Prop $4|q|$ surgeries on $1A$ knots yield their corresponding prism manifolds.
\begin{proof}
Note that the $1A^-$ knots $(\s_1...\s_{2s+2})^{4s+8}(\s_1)^2$ are the twisted torus knots $K(4s+8, 2s+3, 2, 1, 1)$ described in \cite[Theorem 4.1]{JD}. By \cite[Theorem 4.3]{JD}, surface slope surgeries on $K$ yields a Seifert-fibered manifolds with three exceptional fibers of multiplicities $(2, 2, 2s+5)$. Hence by Lemma 2.3 and Remark 2.4, the resulting manifold is $\mcl{P}(2s+5, -2s^2-7s-7)$. Similarly, the $1A^+$ knot $(\s_1...\s_{2s})^{4s}(\s_1)^2$ is $K(4s, 2s+1, 2, 1, 1)$, so the resulting manifold is $\mcl{P}(2s-1, 2s^2+s+1)$.
\end{proof}

\section{Primitive/Seifert-Fibered Knots}
\Def An element $z$ in the free group on two generators $x$ and $y$ is primitive if it is part of a basis, and $(m, n)$ Seifert-fibered if $\langle x, y|z\rangle\cong \langle x, y|x^my^n\rangle$. 

\Def A knot $K$ is primitive  (or resp. $(m, n)$ Seifert-fibered) in $H$ if $w$ is primitive (or resp. $(m, n)$ Seifert-fibered) in $\pi_1(H)$. The property with respect to $H'$ is defined analogously.

\Prop{Knots in families $2, 4, 5$ are (2, 2) Seifert-fibered in $H$ and primitive in $H'$}.
\begin{proof}
Note that for knots $K$ in families $2, 4, 5$, we have $a=1$. This means $w'=x'^my'$ for some $m\in \Z^+$. Since $x'^my'$ and $y'$ form a basis for $\pi_1(H')$, $w'$ is primitive, so $K$ is primitive in $H'$. Now we compute $w$ explicitly for each subfamily of knots.\\

\ti{Family} $5^+_1:$ In this case, $K=(\s_1...\s_{4s+5})^{4ts+6t+4}(\s_1...\s_3)$. Note that $4ts+6t+4=4$ (mod $4s+6$). This means if we start from the $i^{th}$ strand, after passing $x$ once, we will be at the $(i-4)^{th}$ strand (mod $4s+6$). Similarly, traveling along $y$ once sends us from the $i^{th}$ strand to the $(i-1)^{th}$ strand (cf. \autoref{fig:fig2}). Therefore have $w=(x^{s+1}y)^2(x^{s+2}y)^2$ as shown:
\begin{align*}
1 & \xrightarrow{x} 4s+3 \xrightarrow{x} \cdots \xrightarrow{x} 3 \xrightarrow{y} 2 \xrightarrow{x} 4s+4 \xrightarrow{x} \cdots \xrightarrow{x} 4 \xrightarrow{y} 3 \\& \xrightarrow{x} 4s+5  \xrightarrow{x} \cdots \xrightarrow{x} 1 \xrightarrow{y} 4  \xrightarrow{x} 4s+6  \xrightarrow{x} \cdots \xrightarrow{x} 2 \xrightarrow{y} 1
\end{align*}
Here the number $1\leq i\leq 4s+6$ represents a point on the $i^{th}$ strand of the braid. Its position (either before or after $(\s_1...\s_{B-1})^k$) is determined by whether there is a power of $x$ or $y$ behind it. Note that $x^{s+1}y$ and $x^{s+2}y$ form a pair of basis for $\pi_1(H)$. Hence family $5^+_1$ knots are $(2, 2)$ Seifert-fibered with respect to $H$. The computations for other families are similar and we adopt the same notations.\\

\ti{Family} $5^-_1$:\\
$w=(x^{s+1}y)^2(x^{s+2}y)^2$:
\begin{align*}
1 & \xrightarrow{x} 4s+3 \xrightarrow{x} \cdots \xrightarrow{x} 3 \xrightarrow{y} 4 \xrightarrow{x} 4s+6 \xrightarrow{x} \cdots \xrightarrow{x} 2 \xrightarrow{y} 3 \\& \xrightarrow{x} 4s+5  \xrightarrow{x} \cdots \xrightarrow{x} 1 \xrightarrow{y} 2  \xrightarrow{x} 4s+4  \xrightarrow{x} \cdots \xrightarrow{x} 4 \xrightarrow{y} 1
\end{align*}

\ti{Family} $5^-_2$:\\
$w=(x^{s+1}y)^2(x^{s+2}y)^2$:
\begin{align*}
1 & \xrightarrow{x} \cdots  \xrightarrow{x} 4s+5 \xrightarrow{x} 3 \xrightarrow{y} 4 \xrightarrow{x} \cdots \xrightarrow{x}4s+4 \xrightarrow{x} 2 \xrightarrow{y} 3 \\& \xrightarrow{x} \cdots  \xrightarrow{x} 4s+3 \xrightarrow{x} 1 \xrightarrow{y} 2  \xrightarrow{x} \cdots  \xrightarrow{x} 4s+6 \xrightarrow{x} 4 \xrightarrow{y} 1
\end{align*}

\ti{Family} $5^+_2$:\\
$w=(x^{s+2}y)^2(x^{s+3}y)^2$:
\begin{align*}
1 & \xrightarrow{x} \cdots \xrightarrow{x} 4s+9 \xrightarrow{x} 3 \xrightarrow{y} 2 \xrightarrow{x} \cdots \xrightarrow{x} 4s+10 \xrightarrow{x} 4 \xrightarrow{y} 3 \\& \xrightarrow{x} \cdots  \xrightarrow{x} 4s+7 \xrightarrow{x} 1 \xrightarrow{y} 4  \xrightarrow{x} \cdots  \xrightarrow{x} 4s+8 \xrightarrow{x} 2 \xrightarrow{y} 1
\end{align*}

\ti{Family} $4^-_1$:\\
$w=(x^syx^{s+1}yx^{s+1}y)^2(x^{s+1}y)^2$:
\begin{align*}
1 & \xrightarrow{x} 8s-1 \xrightarrow{x} \cdots \xrightarrow{x} 7 \xrightarrow{y} 6 \xrightarrow{x} 8s+4 \xrightarrow{x} \cdots \xrightarrow{x} 4 \xrightarrow{y} 3\\
& \xrightarrow{x} 8s+1 \xrightarrow{x} \cdots \xrightarrow{x} 1 \xrightarrow{y} 8 \xrightarrow{x} 8s+6 \xrightarrow{x} \cdots \xrightarrow{x} 6 \xrightarrow{y} 5\\
& \xrightarrow{x} 8s+3 \xrightarrow{x} \cdots \xrightarrow{x} 3 \xrightarrow{y} 2 \xrightarrow{x} 8s \xrightarrow{x} \cdots \xrightarrow{x} 8 \xrightarrow{y} 7 \\
&  \xrightarrow{x} 8s+5 \xrightarrow{x} \cdots \xrightarrow{x} 5 \xrightarrow{y} 4 \xrightarrow{x} 8s+2 \xrightarrow{x} \cdots \xrightarrow{x} 2 \xrightarrow{y} 1
\end{align*}

\ti{Family} $4^+_1$:\\
$w=(x^{s+1}yx^{s+1}yx^{s+2}y)^2(x^{s+1}y)^2$:
\begin{align*}
1 & \xrightarrow{x} 8s+3 \xrightarrow{x} \cdots \xrightarrow{x} 3 \xrightarrow{y} 4 \xrightarrow{x} 8s+6 \xrightarrow{x} \cdots \xrightarrow{x} 6 \xrightarrow{y} 7\\
& \xrightarrow{x} 8s+9 \xrightarrow{x} \cdots \xrightarrow{x} 1 \xrightarrow{y} 2 \xrightarrow{x} 8s+4 \xrightarrow{x} \cdots \xrightarrow{x} 4 \xrightarrow{y} 5\\
& \xrightarrow{x} 8s+7 \xrightarrow{x} \cdots \xrightarrow{x} 7 \xrightarrow{y} 8 \xrightarrow{x} 8s +10 \xrightarrow{x} \cdots \xrightarrow{x} 2 \xrightarrow{y} 3 \\
&  \xrightarrow{x} 8s+5 \xrightarrow{x} \cdots \xrightarrow{x} 5 \xrightarrow{y} 6 \xrightarrow{x} 8s+8 \xrightarrow{x} \cdots \xrightarrow{x} 8 \xrightarrow{y} 1
\end{align*}

\ti{Family} $4^-_2$:\\
$w=(x^{s+3}yx^{s+2}yx^{s+2}y)^2(x^{s+2}y)^2$:
\begin{align*}
1 & \xrightarrow{x} \cdots \xrightarrow{x} 8s+17 \xrightarrow{x} 7 \xrightarrow{y} 6 \xrightarrow{x} \cdots \xrightarrow{x} 8s+14 \xrightarrow{x} 4 \xrightarrow{y} 3\\
& \xrightarrow{x} \cdots \xrightarrow{x} 8s+11 \xrightarrow{x} 1 \xrightarrow{y} 8 \xrightarrow{x} \cdots \xrightarrow{x} 8s+16 \xrightarrow{x} 6 \xrightarrow{y} 5\\
& \xrightarrow{x} \cdots \xrightarrow{x} 8s+13 \xrightarrow{x} 3 \xrightarrow{y} 2 \xrightarrow{x} \cdots \xrightarrow{x} 8s+18 \xrightarrow{x} 8 \xrightarrow{y} 7 \\
&  \xrightarrow{x} \cdots \xrightarrow{x} 8s+15 \xrightarrow{x} 5 \xrightarrow{y} 4 \xrightarrow{x} \cdots \xrightarrow{x} 8s+12 \xrightarrow{x} 2 \xrightarrow{y} 1
\end{align*}

\ti{Family} $4^+_2$:\\
$w=(x^{s+1}yx^{s+2}yx^{s+2}y)^2(x^{s+2}y)^2$:
\begin{align*}
1 & \xrightarrow{x} \cdots \xrightarrow{x} 8s+9 \xrightarrow{x} 3 \xrightarrow{y} 4 \xrightarrow{x} \cdots \xrightarrow{x} 8s+12 \xrightarrow{x} 6 \xrightarrow{y} 7\\
& \xrightarrow{x} \cdots \xrightarrow{x} 8s+7 \xrightarrow{x} 1 \xrightarrow{y} 2 \xrightarrow{x} \cdots \xrightarrow{x} 8s+10 \xrightarrow{x} 4 \xrightarrow{y} 5\\
& \xrightarrow{x} \cdots \xrightarrow{x} 8s+13 \xrightarrow{x} 7 \xrightarrow{y} 8 \xrightarrow{x} \cdots \xrightarrow{x} 8s+8 \xrightarrow{x} 2 \xrightarrow{y} 3 \\
&  \xrightarrow{x} \cdots \xrightarrow{x} 8s+11 \xrightarrow{x} 5 \xrightarrow{y} 6 \xrightarrow{x} \cdots \xrightarrow{x} 8s+14 \xrightarrow{x} 8 \xrightarrow{y} 1
\end{align*}

\ti{Family} $2^-_1$:\\
$w=(x(xy)^{4s+5})^2(xy)^2$:
\begin{align*}
4s+6\xrightarrow{x} 8s+14 \xrightarrow{xy} \cdots \xrightarrow{xy} 4s+5 \xrightarrow{x} 8s+13 \xrightarrow{xy}\cdots \xrightarrow{xy} 4s+6
\end{align*}

\ti{Family} $2^+_1$:\\
$w=(x(xy)^{4s+3})^2(xy)^2$:
\begin{align*}
4s+8\xrightarrow{x} 8s+10 \xrightarrow{xy} \cdots \xrightarrow{xy} 4s+7 \xrightarrow{x} 8s+9 \xrightarrow{xy}\cdots \xrightarrow{xy} 4s+8
\end{align*}

\ti{Family} $2^-_2$:\\
$w=(x(xy)^{4s+3})^2(xy)^2$:
\begin{align*}
4s+2 \xrightarrow{x} 8s+10 \xrightarrow{xy} \cdots \xrightarrow{xy}4s+1 \xrightarrow{x} 8s+9 \xrightarrow{xy} \cdots \xrightarrow{xy}4s+2
\end{align*}

\ti{Family} $2^+_2:$\\
$w=(x(xy)^{4s+5})^2(xy)^2$:
\be
4s+8 \xrightarrow{x} 8s+14 \xrightarrow{xy} \cdots \xrightarrow{xy}4s+7 \xrightarrow{x} 8s+13 \xrightarrow{xy} \cdots \xrightarrow{xy}4s+8
\qedhere
\ee
\end{proof}
\smallskip

\Prop{Knots in families $1B, 3B$ are (2, 2) Seifert-fibered in $H$ and primitive in $H'$, while the $3A$ knots are (2, 2) Seifert-fibered in $H'$ and primitive in $H$. }\\

\begin{proof}
For these knots, we compute both $w$ and $w'$ explicitly. The computation for the subfamily $1B_1^-$ is described in more details below, while the approaches for the other subfamilies are analogous.\\

\ti{Family} $1B^-_1:$\\
$w=(xy(x^2y)^{2s+2})^2(x^2y)^2$ and $w'=x'(x'^3y'x'^3)^{s+2}x'(x'^3y'x'^3)^{s+1}$, primitive with $x'(x'^3y'x'^3)^{s+1}$:

\begin{align*}
4s+7 \xrightarrow[x']{x} 1 \xrightarrow{y} & 2s+4 \xrightarrow[x'^2]{x}6s+12 \xrightarrow[x']{x} 2s+6 \xrightarrow[y']{y} 1 \xrightarrow[x'^2]{x} 4s+9 \xrightarrow[x']{x} 3 \xrightarrow{y} 2s+6\\
& 2s+6 \xrightarrow[x'^2]{x}6s+14 \xrightarrow[x']{x} 2s+8 \xrightarrow[y']{y} 3 \xrightarrow[x'^2]{x} 4s+11 \xrightarrow[x']{x} 5 \xrightarrow{y} 2s+8\\
& \cdots \\
& 4s+6 \xrightarrow[x'^2]{x}8s+14 \xrightarrow[x']{x} 4s+8 \xrightarrow[y']{y} 2s+3 \xrightarrow[x'^2]{x} 6s+11 \xrightarrow[x']{x} 2s+5 \xrightarrow{y} 4s+8 
\end{align*}

\begin{align*}
4s+8 \xrightarrow[x']{x} 2 \xrightarrow{y} & 2s+5 \xrightarrow[x'^2]{x}6s+13 \xrightarrow[x']{x} 2s+7 \xrightarrow[y']{y} 2 \xrightarrow[x'^2]{x} 4s+10 \xrightarrow[x']{x} 4\xrightarrow{y} 2s+7\\
& 2s+7 \xrightarrow[x'^2]{x}6s+15 \xrightarrow[x']{x} 2s+9 \xrightarrow[y']{y} 4 \xrightarrow[x'^2]{x} 4s+12 \xrightarrow[x']{x} 6\xrightarrow{y} 2s+9\\
& \cdots \\
& 4s+5 \xrightarrow[x'^2]{x}8s+13 \xrightarrow[x']{x} 4s+7 \xrightarrow[y']{y} 2s+2 \xrightarrow[x'^2]{x} 6s+10 \xrightarrow[x']{x} 2s+4 \xrightarrow{y} 4s+7
\end{align*}

\ti{Family} $1B^+_1:$\\
$w=(x^3y(x^2y)^{2s})^2(x^2y)^2$ and $w'=x'^5y'(x'^6y')^sx'^5y'(x'^6y')^{s+1}$, primitive with $x'^5y'(x'^6y')^s$:
\begin{align*}
1 \xrightarrow[x'^5]{x^3}  4s+3 \xrightarrow[y']{y} & 2s+2 \xrightarrow[x'^2]{x} 6s+6 \xrightarrow[x']{x} 2s \xrightarrow{y} 4s+3 \xrightarrow[x'^2]{x}8s+7 \xrightarrow[x']{x} 4s+1 \xrightarrow[y']{y} 2s \\
& \cdots   \\
2 \xrightarrow[x'^5]{x^3} 4s+4 \xrightarrow[y']{y} & 2s+3 \xrightarrow[x'^2]{x} 6s+7 \xrightarrow[x']{x} 2s+1 \xrightarrow{y} 4s+4 \xrightarrow[x'^2]{x}8s+8 \xrightarrow[x']{x} 4s+2 \xrightarrow[y']{y} 2s+1 \\
& \cdots 
\end{align*}

\ti{Family} $3B^-:$\\
$w=(x^2y(xy)^{2s+2})^2(xy)^2$ and $w'=(x'^{2t+1}y'(x'^{2t+2}y')^{s+1})^2x'^{2t+2}y'$, primitive with $x'^{2t+1}y'(x'^{2t+2}y')^{s+1}$:
\begin{align*}
4s+7 \xrightarrow[x'^{2t+1}]{x^2} 1\xrightarrow[y']{y} & 2s+4 \xrightarrow[x'^{t+1}]{x} 2s+6 \xrightarrow{y} 1 \xrightarrow[x'^{t+1}]{x} 3 \xrightarrow[y']{y} 2s+6\\
& \cdots \\
4s+8 \xrightarrow[x'^{2t+1}]{x^2} 2\xrightarrow[y']{y} & 2s+5 \xrightarrow[x'^{t+1}]{x} 2s+7 \xrightarrow{y} 2 \xrightarrow[x'^{t+1}]{x} 4 \xrightarrow[y']{y} 2s+7\\
& \cdots
\end{align*}

\ti{Family} $3B^+:$\\
$w=(x^2y(xy)^{2s})^2(xy)^2$ and $w'=(x'^{2t+3}y'(x'^{2t+2}y')^{s})^2x'^{2t+2}y'$, primitive with $x'^{2t+3}y'(x'^{2t+2}y')^{s}$:
\begin{align*}
1 \xrightarrow[x'^{2t+3}]{x^2}  4s+3 \xrightarrow[y']{y} & 2s+2 \xrightarrow[x'^{t+1}]{x} 2s \xrightarrow{y} 4s+3 \xrightarrow[x'^{t+1}]{x} 4s+1 \xrightarrow[y']{y} 2s\\
& \cdots \\
2 \xrightarrow[x'^{2t+3}]{x^2} 4s+4 \xrightarrow[y']{y} & 2s+3 \xrightarrow[x'^{t+1}]{x} 2s+1 \xrightarrow{y} 4s+4 \xrightarrow[x'^{t+1}]{x} 4s+2 \xrightarrow[y']{y} 2s+1\\
& \cdots
\end{align*}

\ti{Family} $3A^-:$\\
$w'=(x'y'x'^{2s+2})^2x'^2$ and $w=(x^{t+1}y(x^{t+1}yx^ty)^{s+1})^2x^{t+1}yx^ty$, primitive with $x^{t+1}y(x^{t+1}yx^ty)^{s+1}$:
\begin{align*}
2s+4 \xrightarrow[x']{x^{t+1}} 4s+7 \xrightarrow[y']{y} & 1 \xrightarrow[x']{x^{t+1}} 2s+4 \xrightarrow{y} 2s+6 \xrightarrow[x']{x^{t}} 1 \xrightarrow{y} 3 \\    
& \cdots \\
2s+5 \xrightarrow[x']{x^{t+1}} 4s+8 \xrightarrow[y']{y} & 2 \xrightarrow[x']{x^{t+1}} 2s+5 \xrightarrow{y} 2s+7 \xrightarrow[x']{x^{t}} 2 \xrightarrow{y} 4 \\    
& \cdots
\end{align*}

\ti{Family} $3A^+:$\\
$w'=(x'y'x'^{2s+2})^2x'^2$ and $w=(x^{t+1}y(x^{t+1}yx^{t+2}y)^{s+1})^2x^{t+1}yx^{t+2}y$, primitive with $x^{t+1}y(x^{t+1}yx^{t+2}y)^{s+1}$:
\begin{align*}
2s+4  \xrightarrow[x']{x^{t+1}} 1  \xrightarrow[y']{y} & 4s+7 \xrightarrow[x']{x^{t+1}} 2s+4 \xrightarrow{y} 2s+2 \xrightarrow[x']{x^{t+2}} 4s+7 \xrightarrow{y} 4s+ 5 \\
& \cdots \\
2s+5  \xrightarrow[x']{x^{t+1}} 2  \xrightarrow[y']{y} & 4s+8 \xrightarrow[x']{x^{t+1}} 2s+5 \xrightarrow{y} 2s+3 \xrightarrow[x']{x^{t+2}} 4s+8 \xrightarrow{y} 4s+6 \\
& \cdots
\end{align*}

\begin{figure}[H]
\centering
\includegraphics[scale=0.4]{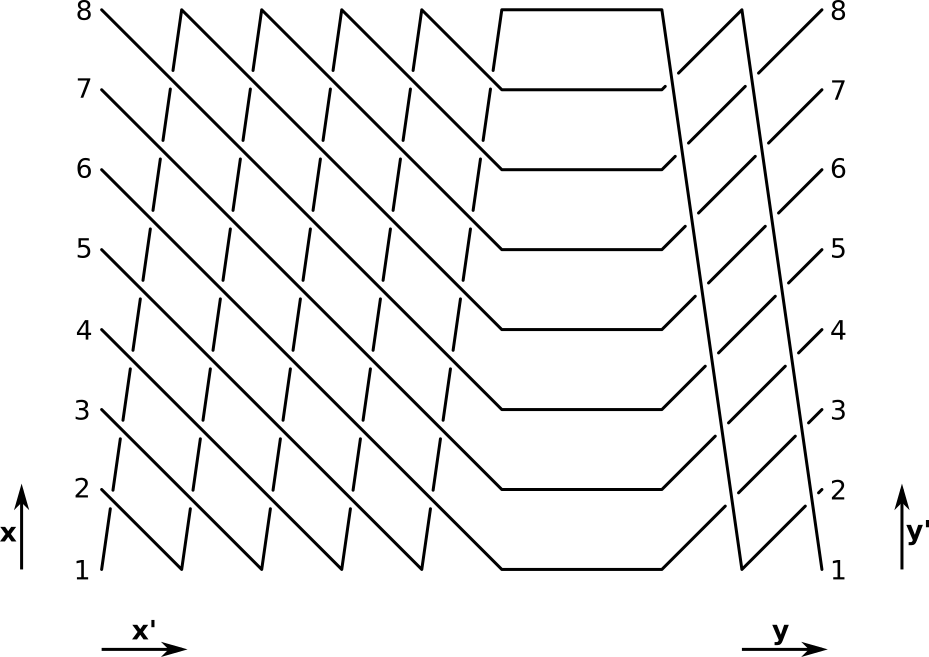}
\caption{The knot $(\s_1...\s_7)^5(\s_7...\s_1)^2$ of family $3A^-$}
\label{fig:fig3}
\end{figure}

\ti{Family} $1B^-_2:$\\
$w=(x^2y(x^3y)^{2s+2})^2(x^3y)^2$ and $w'=x'(x'^2y'x'^2)^{s+1}x'(x'^2y'x'^2)^{s+2}$, primitive with $x'(x'^2y'x'^2)^{s+1}$:
\begin{align*}
1 \xrightarrow[x']{x^2} 4s+7 \xrightarrow{y} & 2s+4 \xrightarrow[x']{x^2} 6s+10 \xrightarrow[x']{x} 2s+2 \xrightarrow[y']{y} 4s+7 \xrightarrow[x']{x^2} 8s+13 \xrightarrow[x']{x} 4s+5 \xrightarrow{y} 2s+2  \\ 
& \cdots   \\
2 \xrightarrow[x']{x^2} 4s+8 \xrightarrow{y} & 2s+5 \xrightarrow[x']{x^2} 6s+11 \xrightarrow[x']{x} 2s+3 \xrightarrow[y']{y} 4s+8 \xrightarrow[x']{x^2} 8s+14 \xrightarrow[x']{x} 4s+6 \xrightarrow{y} 2s+3  \\ 
& \cdots  
\end{align*}

\ti{Family} $1B^+_2:$\\
$w=(x^4y(x^3y)^{2s})^2(x^3y)^2$ and $w'=x'^3y'(x'^4y')^{s+1}x'^3y'(x'^4y')^s$, primitive with $x'^3y'(x'^4y')^s$:
\begin{align*}
4s+3  \xrightarrow[x'^3]{x^4}  1 \xrightarrow[y']{y} & 2s+2 \xrightarrow[x']{x^2} 6s+8 \xrightarrow[x']{x} 2s+4 \xrightarrow{y} 1 \xrightarrow[x']{x^2} 4s+7 \xrightarrow[x']{x} 3 \xrightarrow[y']{y} 2s+4 \\
& \cdots \\
4s+4 \xrightarrow[x'^3]{x^4} 2 \xrightarrow[y']{y} & 2s+3 \xrightarrow[x']{x^2} 6s+9 \xrightarrow[x']{x} 2s+5 \xrightarrow{y} 2 \xrightarrow[x']{x^2} 4s+8 \xrightarrow[x']{x} 4 \xrightarrow[y']{y} 2s+5\\
& \cdots \qedhere
\end{align*} 
\end{proof}

\Prop{$4|q|$ surgeries on knots in families $1B, 2, 3A, 3B, 4$, and 5 yield prism manifolds.}

\begin{proof}
By Propositions 3.4 and 3.5, every such knot $K$ is primitive/(2, 2)-Seifert-fibered. Thus \cite[Proposition 2.3]{JD} implies surface slope surgery on $K$ yields a manifold $\mcl{P}$ that is either Seifert-fibered over $S^2$ with three exceptional fibers or a connected sum of two lens spaces. If $\mcl{P}$ is a connected sum of two lens spaces, since the two exceptional fibers have multiplicities 2 and 2, $\mcl{P}=L(2, 1) + L(2, 1)$. Then $|H_1(\mcl{P})|=4$, which means the surgery coefficient is 4. However, $|q|\neq 1$ for each family in Table 1. Contradiction! Consequently, $\mcl{P}$ is a Seifert-fibered space over $S^2$ with three exceptional fibers. Since the first two multiplicities are 2 and 2, $\mcl{P}$ is a prism manifold.
\end{proof}

\section{Multiplicity of The Third Exceptional Fiber}
In this section we compute the absolute value of the multiplicity of the third exceptional fiber of $\mcl{P}$ resulting from surgeries on knots in families $1B, 2, 3A, 3B, 4$, and 5. By \cite[Lemma 6.2]{NZ}, we have:

\Lem If $\mcl{P}(p, q)$ is realizable by $4|q|$ surgery on a knot $K$, then $|p|=|\Delta_K(-1)|$.\\

Hence determining $|p|$ amounts to computing the determinant of $K$. Note that for all knots $K$ in the families above, we have odd $B$. Since $K$ and $K\s_{B+1}$ represent the same knot, it is equivalent to compute $|\det(K\s_{B+1})|$. Observe that the braid word $K\s_{B+1}$ has $B+2$ strands. Consider the reduced Burau representation $\ovl{\mfk{b}}$ of the braid group on $B+2$ strands given by 

\be
\ovl{\mfk{b}}(\s_1)=
\begin{pmatrix}
-\mfk{t} & 1 & 0 \\
0 & 1 & 0 \\
0 & 0 & I_{B-1}
\end{pmatrix}
\ \ \ \
\ovl{\mfk{b}}(\s_{i})=
\begin{pmatrix}
I_{i-2} & 0 & 0 & 0 & 0\\
0 & 1 & 0 & 0 & 0 \\
0 & \mfk{t} & -\mfk{t}  & 1 & 0 \\
0 & 0 & 0 & 1 & 0\\
0 & 0 & 0 & 0 & I_{B-i}
\end{pmatrix}
\ \ \ \
\ovl{\mfk{b}}(\s_{B+1})=
\begin{pmatrix}
I_{B-1} & 0 & 0 \\
0 & 1 & 0 \\
0 & \mfk{t} & -\mfk{t}
\end{pmatrix}
\ee

where $2\leq i\leq B$ and $I_m$ is the $m\times m$ identity matrix. Note that we have 
\be
\Delta_{K\s_{B+1}}(\mfk{t})=\frac{1-\mfk{t}}{1-\mfk{t}^{B+2}}\det(I-\ovl{\mfk{b}}(K\s_{B+1}))\ee 
Let $\mfk{b}$ be the restriction of $\ovl{\mfk{b}}$ to $\mfk{t}=-1$. Then since $B+2$ is odd, 
\be
|\Delta_K(-1)|=|\Delta_{K\s_{B+1}}(-1)|=|\det(I-\mfk{b}(K\s_{B+1}))|
\ee
This provides a straightforward approach to compute $|\Delta_K(-1)|$. For matrices in this section, blank means all zeros, while horizontal, vertical, or diagonal dots refer to patterns along rows, columns, or diagonals. (For diagonal dots, all non diagonal entries in the block are zeros). In terms of patterns, $(m, m, \cdots)$ refers to a line of $m$'s, while $(m, n, \cdots)$ means a line of alternating $m$'s and $n$'s. (Here vertical dots and diagonal dots go from top to bottom.)
We first make a few simple observations about the representations of half twists and full twists of $\s_1...\s_B$.\\

\Lem{$\mfk{b}(\s_1...\s_B)$ is the $(B+1)\times (B+1)$ matrix}

\be
\begin{pmatrix}
0 & 0 &  & 0 & 1 & 1 \\
-1 & 0 &   & 0 & 1 & 1 \\
0 & -1 &   & 0 & 1 &1\\
 &  & \ddots &  & \vdots & \vdots \\
0 & 0 &  & -1 & 1 & 1\\
0 & 0 &  & 0 & 0 & 1
\end{pmatrix}
\ee
\smallskip

\begin{proof}
For $2\leq i\leq B+1$, let $M_i$ be the $i\times i$ lower right part of the matrix above. Then by induction, 
\be
\mfk{b}(\s_1...\s_i)=\begin{pmatrix} M_{i+1} & 0 \\ 0 & I_{B-i}
\end{pmatrix}
\ee
for $1\leq i\leq B$. Hence in particular, $\mfk{b}(\s_1...\s_B)$ is the matrix specified above.
\end{proof}

\Lem $\mfk{b}((\s_1...\s_B)^{m\frac{B+1}{2}})$ is the matrix on the left when $m$ is even and on the right when $m$ is odd. (The quadrants of the matrix on the right correspond to the quadrants of the actual matrix.)
\be
\begin{pmatrix}
1 & 0  &  & 0 & m \\
0 & 1  &  & 0 & 0 \\
 &  &  \ddots &  & \vdots \\
0 & 0   &  & 1 & m \\
0 & 0   & & 0 & 1 \\
\end{pmatrix}
\ \ \ \ \ \ \ \ 
\begin{pmatrix}
0 & 0 &  & 0 & 1 & -1 & 0 & & 0 & m+1 \\
0 & 0 &   & 0 &  1& 0 & -1 & & 0 & 0 \\

 & &  & & \vdots& & & \ddots  &  & \vdots\\

0 & 0 & & 0 & 1 & 0 & 0 & & -1 & 0 \\
0 & 0 & & 0 & 1 & 0 & 0 & & 0 & m \\
-1 & 0 & & 0 & 1 & 0 & 0 & & 0 & 1 \\
0 & -1 & & 0 &  1 & 0 & 0 & & 0 &  m\\

 & &  \ddots  & & \vdots & & & & & \vdots\\

0 & 0 & & -1 & 1 & 0 & 0 & & 0 & m \\
0 & 0 & & 0 & 0 & 0 & 0 & & 0 & 1 \\
\end{pmatrix}
\ee
\smallskip

\Rem Lemma 4.3 is proved by straightforward inductions, which we omit for concision. Note that for all knots $K$ in families $1B, 2, 3A, 3B, 4$, and $5$, we have $k$ congruent to some constants (mod $\frac{B+1}{2}$). Hence this lemma greatly simplifies the computation of $\mfk{b}(K\s_{B+1})$.
\smallskip

\Prop For $K$ in families $1B, 2, 3A, 3B, 4, 5$  and their corresponding $\mcl{P}(p, q)$, we have
\be
|\det(I-\mfk{b}{K\s_{B+1}})|=|p|
\ee

\newpage
\begin{proof}
We illustrate the proof approach with the $5^+_1$ knots $(\s_1...\s_{4s+5})^{4ts+6t+4}\s_1\s_2\s_3$. Note that 
\be
\mfk{b}((\s_1...\s_{4s+5})^{4ts+6t+4}\s_1\s_2\s_3\s_{4s+6})=\mfk{b}((\s_1...\s_{4s+5})^{2t(2s+3)})\mfk{b}((\s_1...\s_{4s+5})^4)\mfk{b}(\s_1\s_2\s_3\s_{4s+6})
\ee
We can compute $\mfk{b}(\s_1\s_2\s_3\s_{4s+6})$ (left) and $\mfk{b}((\s_1...\s_{4s+5})^4)$ (right) directly:
\be
\begin{pmatrix}
0 & 0 & 1 & 1 & & 0 & 0 \\
-1 & 0 & 1 & 1 & & 0 &0 \\
0 & -1 & 1 & 1 & & 0 &0 \\
0 & 0 & 0 &1 & & 0 & 0  \\
 &  & &  & \ddots & & \\
0 & 0 & 0 & 0  & & 1 & 0 \\
0 & 0 & 0 & 0 & & -1 & 1  \\
\end{pmatrix}
\ \ \ \ \ \ \ \ 
\begin{pmatrix}
0 & 0 &  & 0& -1 & 1 & 0 & 0 & 2 \\
0 & 0 &   & 0& -1 & 0 & 1 & 0 &0 \\
0 & 0 &   & 0& -1 & 0 & 0  & 1 & 2 \\
0 & 0 &  & 0& -1 & 0 & 0 & 0 & 1  \\
1 &0 &   & 0& -1 & 0  & 0& 0 & 1 \\
0 &1 &   & 0& -1 & 0  & 0& 0 & 1 \\
& & \ddots & & \vdots & & & & \vdots  \\
0 &0 &  &1 & -1 & 0 &0  & 0 & 1  \\
0 &0 &  & 0& 0 & 0 &0  & 0 & 1  \\
\end{pmatrix}
\ee
Then by Lemma 4.3 , we have $I_{4s+6}-\mfk{b}(K\s_{B+1})=$
\be
\begin{pmatrix}
1 & & & & & & & &  & 0 & 1 & -1 & 0 & 2t+2 & -2t-2 \\
& 1& & & & & & & & 0 & 1 & 0 & -1 & 0 & 0 \\
& & 1& & & & & &  & 0 & 1 & 0 & 0 &2t+1 & -2t-2 \\
& & & 1& & & & & & 0 & 1 & 0 & 0 & 1 & -1\\
0 & 0 & -1& -1& \ddots& & & &  & 0 & 1 & 0& 0& 2t+1& -2t-1\\
1& 0 &-1 &-1 & & \ddots& & &  & 0 & 1 & 0& 0& 1& -1\\
0 & 1& -1& -1& & & \ddots& & &  & \vdots & & & \vdots & \vdots\\
0 & 0 & 0 & -1& & & & \ddots& &  & \vdots & & & \vdots& \vdots\\

& & & &\ddots & & & & \ddots & & \vdots& & & \vdots& \vdots\\

& & & & & \ddots & & & &1 & 1 & 0 & 0 & 2t+1 & -2t-1 \\
& & & & & &\ddots & & & 0& 2 & 0 & 0 & 1 & -1 \\
& & & & & & & \ddots &  & 0 & 1 & 1& 0& 2t+1 & -2t-1 \\
& & & & & & & &  \ddots & 0 & 1 & 0 & 1 & 1 & -1 \\
& & & & & & & &  & -1 & 1 & 0 & 0 & 2t+2 & -2t-1 \\
& & & & & & & &  & 0 & 0 & 0 & 0 & 1 & 0\\
\end{pmatrix}
\ee 
Adding the last column to the second to last one and  reducing every four rows, we obtain an upper triangular block matrix with $I_{4s}$ and the matrix below in the upper left and lower right corners.
\be
\begin{pmatrix}
1 & s+2 & 0 & 0 & -1 & -2ts-2t -s-3\\
0 & s+2 & 1 & 0 & -1 & -s-1\\
0 & s+1 &1 & 1 & -1 & -2ts - 2t -s -3\\
0 & s+1 & 0& 1 & 0 & -s-1 \\
-1 & 1 & 0 & 0 & 1 & -2t-1 \\
0 & 0 & 0 & 0 & 1 & 0 \\
\end{pmatrix}
\ee
Then it is easy to verify that  $|\det(I-\mfk{b}(K\s_{B+1}))=4s^2t+16st+14t+4s+11$ for family $5^+_1$.
\end{proof}

\Rem Our proof of Proposition 4.5 for other families is identical to the case above. This proposition, combined with Lemma 2.3 and Proposition 3.5, completes the proof of Theorem 1.1.

\noindent\\
Department of Mathematics, California Institute of Technology, Pasadena, CA 91125\\
E-mail address: zshang@caltech.edu


\begin{thebibliography}{9}

\bibitem{B1}
{W. Ballinger, C. C. Hsu, W. Mackey, Y. Ni, T. Ochse, F. Vafaee. \ti{The
prism manifold realization problem}, Preprint, arXiv:1612.04921, 2016.}

\bibitem{B2}
{W. Ballinger, Y. Ni, T. Ochse, F. Vafaee. \ti{The prism manifold realization
problem II}, Preprint, arXiv:1710.00089, 2017.}

\bibitem{B3}
{W. Ballinger, Y. Ni, T. Ochse, F. Vafaee. \ti{The prism manifold realization
problem III}, Preprint, arXiv:1808.05321, 2018.}

\bibitem{BH}
{S. Bleiler, C. Hodgson. \ti{Spherical space forms and Dehn filling}, Topology, 35(3):809-833, 1996.}

\bibitem{BL}
{S. Boyer, D. Lines. \ti{Surgery formulae for Casson's invariant and extensions to homology lens spaces}, J. Reine Angew. Math., 405:181-220, 1990. }

\bibitem{BZ}
{S. Boyer, X. Zhang. \ti{Finite Dehn surgery on knots}, J. Amer. Math. Soc., 9:1005-1050, 1996.}


\bibitem{JD}
{J. C. Dean. \ti{Small Seifert-fibered Dehn surgery on hyperbolic knots}, Algebr. Geom. Topol., 3:435-472, 2003.}

\bibitem{JG}
{J. E. Greene. \ti{The lens space realization problem}, Ann. Math., 177(2): 449-511, 2013.}

\bibitem{LG}
{L. Gu. \ti{Integral finite surgeries on knots in $S^3$}, Preprint, arXiv:1401.6708, 2014.}

\bibitem{CL}
{C. Lescop. \ti{Global surgery formula for the Casson-Walker invariant}, (AM-140), Princeton University Press, 1996.}

\bibitem{LN}
{E. Li, Y. Ni. \ti{Half-intergal finite surgeries on knots in $S^3$}, Ann. Fac. Sci. Toulouse Math., (6), 24(5):1157-1178, 2015. }

\bibitem{WL}
{W. B. R. Lickorish. \ti{A representation of orientable combinatorial 3-manifolds}, Ann. Math., 76:531-540, 1962.}

\bibitem{LM}
{L. Moser. \ti{Elementary surgery along a torus knot}, Pac. J. Math., 38:737-745, 1971.}

\bibitem{NZ}
{Y. Ni, X. Zhang. \ti{Finite Dehn surgeries on knots in $S^3$}, Algebr. Geom. Topol., 18:441-492, 2018.}

\bibitem{AW}
{A. H. Wallace. \ti{Modifications and cobounding manifolds}, Canad. J. Math., 12:503-528, 1960.}
\end{thebibliography}
\end{document}